\documentclass[8 pt]{amsproc}
\usepackage{amssymb}
\oddsidemargin=-1cm\evensidemargin=-1cm
\begin{document}
\numberwithin{equation}{section}

\def\1#1{\overline{#1}}
\def\2#1{\widetilde{#1}}
\def\3#1{\widehat{#1}}
\def\4#1{\mathbb{#1}}
\def\5#1{\frak{#1}}
\def\6#1{{\mathcal{#1}}}

\def\C{{\4C}}
\def\R{{\4R}}
\def\N{{\4N}}
\def\Z{{\4Z}}

\title{Normal Forms and Degenerate CR Singularities}  

\author{  Valentin Burcea  }

\begin{abstract}  We construct a normal form for a class of real-formal surfaces $M\subset\mathbb{C}^{2}$ defined near a degenerate CR singularity $p=0$ by
$$w=P\left(z,\overline{z}\right) + \mbox{O}\left(\left|z\right|^{k_{0}+1}\right),$$
 where $P\left(z,\overline{z}\right)$  is a real-valued  homogeneous polynomial in  $\left(z,\overline{z}\right)$ of degree $k_{0} \geq 3 $ such that the coefficients of $z^{k_{0}}$ and  $\overline{z}^{k_{0}}$
 are vanishing. 
\end{abstract}
\address{V. Burcea: Departament of Mathematics, The University Federal of Santa Catarina, Florianopolis, Brazil}

\email{vdburcea@gmail.com}


\thanks{This
project has been supported partially by CAPES}

\maketitle

\def\cn{{\C^n}}
\def\cnn{{\C^{n'}}}
\def\ocn{\2{\C^n}}
\def\ocnn{\2{\C^{n'}}}

% Abbreviations

\def\dist{{\rm dist}}
\def\const{{\rm const}}
\def\rk{{\rm rank\,}}
\def\id{{\sf id}}
\def\tr{{\bf tr\,}}
\def\aut{{\sf aut}}
\def\Aut{{\sf Aut}}
\def\CR{{\rm CR}}
\def\GL{{\sf GL}}
\def\Re{{\sf Re}\,}
\def\Im{{\sf Im}\,}
\def\span{\text{\rm span}}
\def\Diff{{\sf Diff}}

\def\codim{{\rm codim}}
\def\crd{\dim_{{\rm CR}}}
\def\crc{{\rm codim_{CR}}}

\def\phi{\varphi}
\def\eps{\varepsilon}
\def\d{\partial}
\def\a{\alpha}
\def\b{\beta}
\def\g{\gamma}
\def\G{\Gamma}
\def\D{\Delta}
\def\Om{\Omega}
\def\k{\kappa}
\def\l{\lambda}
\def\L{\Lambda}
\def\z{{\bar z}}
\def\w{{\bar w}}
\def\Z{{\1Z}}
\def\t{\tau}
\def\th{\theta}

\emergencystretch15pt \frenchspacing

\newtheorem{Thm}{Theorem}[section]
\newtheorem{Cor}[Thm]{Corollary}
\newtheorem{Pro}[Thm]{Proposition}
\newtheorem{Lem}[Thm]{Lemma}

\theoremstyle{definition}\newtheorem{Def}[Thm]{Definition}

\theoremstyle{remark}
\newtheorem{Rem}[Thm]{Remark}
\newtheorem{Exa}[Thm]{Example}
\newtheorem{Exs}[Thm]{Examples}

\def\bl{\begin{Lem}}
\def\el{\end{Lem}}
\def\bp{\begin{Pro}}
\def\ep{\end{Pro}}
\def\bt{\begin{Thm}}
\def\et{\end{Thm}}
\def\bc{\begin{Cor}}
\def\ec{\end{Cor}}
\def\bd{\begin{Def}}
\def\ed{\end{Def}}
\def\br{\begin{Rem}}
\def\er{\end{Rem}}
\def\be{\begin{Exa}}
\def\ee{\end{Exa}}
\def\bpf{\begin{proof}}
\def\epf{\end{proof}}
\def\ben{\begin{enumerate}}
\def\een{\end{enumerate}}
\def\beq{\begin{equation}}
\def\eeq{\end{equation}}

 \section{Introduction and Main Result}
 The study of  real submanifolds in  complex space
near a CR singularity  goes back  to    Bishop\cite{Bi}.  A point $p\in M$  is called a CR
singularity if it is a   discontinuity point for the map $M\ni q\longrightarrow  \dim_{\mathbb{R}}T^{c}_{q}M$  defined near $p$.  Bishop\cite{Bi} considered the case when there exist  coordinates $(z,w)$ in $\mathbb{C}^{2}$ such that near the CR
singularity $p=0$, the surface $M\subset\mathbb{C}^{2}$ is defined  by
\begin{equation}
 w=z\overline{z}+\lambda\left(z^{2}+\overline{z}^{2}\right)+\rm{O}(3),\label{clasic}
 \end{equation}
where  $\lambda\in\left[0,\infty\right]$ is a holomorphic
invariant called the Bishop invariant. When $\lambda=\infty$, $M$
is understood to be defined by the equation
$w=z^{2}+\overline{z}^{2}+\rm{O}(3)$. If $\lambda$ is
non-exceptional,  Moser-Webster\cite{MW} proved that there
exists a formal transformation that sends $M$ into the following normal form
\begin{equation*}
w=z\overline{z}+\left(\lambda+\epsilon
u^{q}\right)\left(z^{2}+\overline{z}^{2}\right),\quad
\epsilon\in\left\{0,-1,+1\right\},\quad q\in\mathbb{N},
\end{equation*}
where $w=u+iv$.   Moser\cite{Mos} constructed when $\lambda=0$
the following partial normal form:
\begin{equation}
 w=z\overline{z}+2\Re\left\{\displaystyle\sum_{j\geq s}a_{j}z^{j}\right\}.\label{moser}
\end{equation}

Here $s:=\min\left\{j\in\mathbb{N}^{\star};\hspace{0.1
cm}a_{j}\neq 0\right\}$ is the simplest higher order invariant,
known as the Moser invariant.   When
$s<\infty$, Huang-Yin\cite{HY2}  proved that
\eqref{moser} can be formally transformed into the following normal form 
\begin{equation}
w=z\overline{z}+ 2\Re\left\{\displaystyle\sum_{j\geq
s}a_{j}z^{j}\right\},\quad a_{s}=1,\quad a_{j}=0,\quad
\mbox{if}\quad j=0,1\hspace{0.1 cm}\mbox{mod}\hspace{0.1
cm}s,\quad j>s.\label{hye}
\end{equation}

In this note, we construct a normal form for a  surface $M\subset\mathbb{C}^{2}$
 defined near $p=0$ by
 \begin{equation}w=P(z,\overline{z})+\displaystyle\sum_{j+l\geq k_{0}+1}a_{j,l}z^{j}\overline{z}^{l},\label{wer}\end{equation}
where $P(z,\overline{z})$ is  real-valued  homogeneous polynomial in  $\left(z,\overline{z}\right)$ of degree $k_{0} \geq 3 $ having the coefficients of its pure terms vanishing. 

Our case (\ref{wer}) is different from the classical case from
(\ref{clasic}) when $\lambda=0$, studied intially by Moser\cite{Mos}, and requires a different approach using the Fischer decomposition\cite{Sh}, that has been applied by Zaitsev\cite{D1},\cite{D2},\cite{D3} in other situations. In order  to
develop the partial normal form,  we introduce the  space 
\begin{equation}S_{N},\label{space} \end{equation}
which consists homogeneous polynomials of degree $N$ in $(z,\overline{z})$, denoted as 
$S(z,\overline{z})=S_{0}(z,\overline{z})$,   such that
$$S_{j}(z,\overline{z})=S_{j+1}(z,\overline{z})P(z,\overline{z})+T_{j+1}(z,\overline{z}),\quad P^{\star}\left(T_{j+1}(z,\overline{z})\right)=0,$$
 where the following normalizations hold:
 $$\left(z^{N-jk_{0}}\right)^{\star}\left(S_{j}(z,\overline{z})\right)=0,\hspace{0. cm}\forall j=0,\dots,\left[\frac{N}{k_{0}}\right],\quad  \left(\overline{z}^{N-jk_{0}}\right)^{\star}\left(S_{j}(z,\overline{z})\right)=0,\hspace{0.1 cm}\forall\hspace{0.1 cm}j=1,\dots,\left[\frac{N}{k_{0}}\right],$$
for all $N\geq k_{0}+1$, where we have used by \cite{Sh} the following notation
\begin{equation}P^{\star}=\displaystyle\sum_{m+n=k_{0}}\overline{p}_{m,n}\frac{\partial^{m+n}}{\partial z^{m}\partial \overline{z}^{n}},\quad\mbox{provided $P(z,\overline{z})=\displaystyle\sum_{m+n=k_{0}}p_{m,n}z^{m}\overline{z}^{n}$.}\label{pol}\end{equation}

We observe that the Fischer normalization conditions defined by $\mathcal{S}_{N}$ generalize the normalization conditions defined by the Moser partial normal form \cite{Mos}. More exactly, the space $\mathcal{S}_{N}$ defined by (\ref{space}), with respect to the polynomial $P(z,\overline{z})$, is trivial in Moser's case\cite{Mos}.

In order to find a normal form for the   surfaces defined by (\ref{wer}), we impose firstly the following nondegeneracy condition
\begin{equation}s:=\min \left\{a_{0,l};\quad l\geq k_{0}+1\right\}<\infty.\label{nondeg}\end{equation}

Here $s$  is the simplest holomorphic invariant, similar to the classical case from (\ref{clasic}) when $\lambda=0$, which was studied by Huang-Yin\cite{HY2}. The nondegeneracy condition (\ref{nondeg}) helps us to find a normal form for the   surfaces defined as in (\ref{wer}). Our case requires also the following nondegeneracy conditions
\begin{equation}0\neq p_{1,k_{0}-1},\quad  \alpha \neq 0,s,\quad\mbox{where $zP_{z}(z,\overline{z})=\alpha P(z,\overline{z})+R(z,\overline{z})$ and  $P^{\star}\left(R(z,\overline{z})\right)=0$. }\label{alpha}\end{equation} 

Throughout this note, in order to simplify the computations, we assume
 $p_{1,k_{0}-1}=1$. The main result  is the following

 \bt \label{teo}Let $M\subset\mathbb{C}^{2}$ be a formal surface defined near  $p=0$ by (\ref{wer}) satisfying the nondegeneracy conditions (\ref{nondeg}) and (\ref{alpha}). Then there exists a unique formal transformation of the following type
\begin{equation}\left(z',w'\right)=\left(z+\displaystyle\sum_{k+l\geq 2}f_{k,l}z^{k}w^{l},\quad w+\displaystyle\sum_{k+l\geq 2}g_{k,l}z^{k}w^{l}\right) ,\label{2800}\end{equation}
that transforms $M$ into the following formal normal form:
\begin{equation}w'=P\left(z',\overline{z'}\right) +\displaystyle\sum _{m+n\geq k_{0}+1 }a'_{m,n}{z'}^{m}\overline{z'}^{n},\label{3000}
\end{equation}
where  the following Fischer
normalization conditions are satisfied
\begin{equation}a'_{N,0}=\overline{a'_{N,0}}\hspace{0.1 cm}\mbox{and}\hspace{0.1 cm}\displaystyle\sum _{m+n=N\atop{m,n\neq 0}}a'_{m,n}{z}^{m}\overline{z}^{n}\in\mathcal{S}_{N},\quad\mbox{for all $N\geq k_{0}+1$,} \label{cn}
\end{equation}
where $\mathcal{S}_{N}$ is defined in (\ref{space}), and as well the following  normalization conditions hold
\begin{equation}a'_{0,k}=0,\quad\mbox{for all $k\equiv 0, k_{0}-1\mod s$.} \label{extra}\end{equation} \et

The normalization conditions  (\ref{cn}) leave  undetermined an infinite number of parameters   acting at the higher degree levels in
the local defining equations, because the group of formal  automorphisms preserving the origin of the model $w=P(z,\overline{z})$ is infinite dimensional
similarly to the classical case from (\ref{clasic}) when $\lambda=0$ studied by Huang-Yin\cite{HY2}. The space $\mathcal{S}_{N}$ defined by (\ref{space}) helps us just to develop a partial normal form
and those free parameters are determined using the normalization conditions  (\ref{extra}). In order to impose these normalization conditions, we apply the methods developed by Huang-Yin\cite{HY2} using the model $w=P(z,\overline{z})+z^{s}+\overline{z}^{s}$ instead of the previous model   $w=P(z,\overline{z})$. The only obstacle
that appears here is that we can not define a system of weights for $(z, \overline{z})$ as in the case of Huang-Yin\cite{HY2}, that can make the model $w=P(z,\overline{z})+z^{s}$
homogeneous. In order to overcome this problem, we  use a different strategy by considering a different type of weights system called
here system of pseudo-weights (see Subsection $4.1$), which helps us to apply the methods of Huang-Yin\cite{HY2} preserving the normalization
conditions defined by the space $\mathcal{S}_{N}$. These methods allow us to construct a normal form, but can not detect any complete set of simple invariants associated
to the  surface $M$ defined by (\ref{3000}). We   mention that Huang and Yin discovered a complete system of invariants for the surfaces
defined by (\ref{clasic}) with the vanishing Bishop invariant. They\cite{HY2} proved that the set of complex numbers $\left\{a_{i}\right\}_{i\geq s}$ given by (\ref{hye}) defines a complete set of invariants for (\ref{hye}).

We  mention  that   real submanifolds, near a CR singularity in   complex spaces, under unimodular transformations,
have been studied by Gong\cite{G1},\cite{G2}. Furthermore, the CR singular points of the real submanifolds are interesting because the CR singularities
can contribute to the structure of the local hull of holomorphy as it has been shown by Kenig-Webster\cite{KW}, Huang-Krantz\cite{HK} and Huang-Yin
\cite{HYA},\cite{HYB}.   Forstneric-Stout\cite{FS} proved that the   surfaces defined in (\ref{clasic}) with $2\lambda>1$ are local polynomial convex at $p=0$. Problems regarding the local polynomial convexity when  $2\lambda=1$  have been studied by J{\"o}ricke\cite{Jo}. Surfaces in $\mathbb{C}^{2}$ with degenerate CR singularities and related local polynomial convexity problems near degenerate CR singularities have
been studied by Bharali\cite{B1},\cite{B2},\cite{B3},\cite{B4}. In codimensions different from $2$, we mention the work of Coffman\cite{C1},\cite{C2} about CR singularities.

We do not know if the normal form (\ref{3000}) is convergent. In general,   a non-trivial normal form is not  necessarily  convergent  as it has been shown recently by Kolar\cite{K2}. For cases when a normal form is convergent, we  mention the recent work of Kossovskiy-Zaitsev\cite{KZ}.

\section{Preliminaries}
\subsection{ Notations}  Throughout this note, we  use   notations of the following type
$$a_{\geq l}(z,\overline{z})=\displaystyle\sum_{m+n\geq l}a_{m,n}z^{m}\overline{z}^{n}\hspace{0.1 cm}\mbox{and}\hspace{0.1 cm} a_{l}(z,\overline{z})=\displaystyle\sum_{m+n=l}a_{m,n}z^{m}\overline{z}^{n},\quad \mbox{for all $l\geq k_{0}+1$}.$$

\subsection{  Transformation Equations}  Let $M\subset\mathbb{C}^{2}$ be a formal surface defined near
$p=0$ by  
\begin{equation}w=P(z,\overline{z})+\displaystyle\sum_{m+n\geq k_{0}+1}a_{m,n}z^{m}\overline{z}^{n}.\label{a1}\end{equation}

Let $M'\subset\mathbb{C}^{2}$ be another formal  surface defined near $p'=0$ by  
\begin{equation}w'=P\left(z',\overline{z'}\right)+\displaystyle\sum_{m+n\geq k_{0}+1}a'_{m,n}{z'}^{m}\overline{z'}^{n}.\label{a2}\end{equation}

Let $\left(z',w'\right)= \left(f(z,w), g(z,w)\right)$ be a formal equivalence transforming $M$ into $M'$ and  fixing the point $p=0\in\mathbb{C}^{2}$. It
follows by (\ref{a2}) that
\begin{equation}g(z,w)=P\left(f(z,w),\overline{f(z,w)}\right)+\displaystyle\sum_{m+n\geq k_{0}+1}a'_{m,n}\left(f(z,w)\right)^{m}\overline{\left(f(z,w)\right)^{n}},\label{a3}\end{equation}
where $w$ is defined by (\ref{a1}).  Writing  $f(z,w)=\displaystyle\sum_{ m+n\geq 0}
f_{m,n}z^{m}w^{n}$ and $g(z,w) = \displaystyle\sum_{ m+n\geq 0}
g_{m,n}z^{m}w^{n}$,   it follows by (\ref{a3}) that

\begin{equation}\begin{split} &\quad\quad\quad\quad\quad\quad\quad\quad\quad\displaystyle\sum_{ m+n\geq 0}
g_{m,n}z^{m}\left(P(z,\overline{z})+a_{\geq k_{0}+1}(z,\overline{z})\right)^{n}=  P\left(\displaystyle\sum_{ m+n\geq 0}
f_{m,n}z^{m}\left(P(z,\overline{z})+a_{\geq k_{0}+1}(z,\overline{z})\right)^{n},\right.\\&\left.
\quad\quad\quad\quad\quad\quad\quad\quad\quad\quad\quad\quad\quad\quad\quad\quad\quad\quad\quad\quad\quad\quad\quad\quad\quad\quad\quad\quad\quad\quad \overline{\displaystyle\sum_{ m+n\geq 0}
f_{m,n}z^{m}\left(P(z,\overline{z})+a_{\geq k_{0}+1}(z,\overline{z})\right)^{n}}\right)\\&\quad\quad\quad\quad\quad\quad+a'_{\geq k_{0}+1}\left(\displaystyle\sum_{ m+n\geq 0}
f_{m,n}z^{m}\left(P(z,\overline{z})+a_{\geq k_{0}+1}(z,\overline{z})\right)^{n},\overline{\displaystyle\sum_{ m+n\geq 0}
f_{m,n}z^{m}\left(P(z,\overline{z})+a_{\geq k_{0}+1}(z,\overline{z})\right)^{n}}\right).\end{split} \label{b1}\end{equation}

Since our map fixes the point $p=0\in\mathbb{C}^{2}$, it follows that 
 $g_{0,0}=0$ and $f_{0,0}=0.$ Collecting the terms of bidegree $(m, 0)$ in $(z, \overline{z})$ in (\ref{b1}), for all $m<k_{0}$, it follows that  $g_{m,0}= 0$, for all $m<k_{0}$.  Collecting   the sum  of terms of bidegree $(m, n)$ in $(z, \overline{z})$ with $m + n = k_{0}$ in (\ref{b1}), it
follows that
\begin{equation} g_{0,1}P(z, \overline{z}) = P \left(f_{1,0}z, \overline{f_{1,0}z}\right).\label{b2} \end{equation}

Then (\ref{b2}) describes all the possible values of $g_{0,1}$ and $f_{1,0}$. In particular, 
we obtain that $\Im g_{0,1}=0$. By composing with a linear
automorphism of the Model $  w = P(z, \overline{z})$,  we can assume that $g_{0,1}=1$ and $f_{1,0}=1$.   Then, by a careful analysis of the terms interactions in (\ref{b1}), we conclude that in order to put suitable normalization conditions, we have to
consider the following terms
$$ g_{m,n}z^{m}\left(P(z, z)\right)^{n} ,\quad f_{m,n}z^{m}P_{z}(z, \overline{z})\left(P(z, \overline{z})\right)^{n} ,\quad \overline{f_{m,n}z^{m}P_{z}(z, \overline{z})}\left(P(z, \overline{z})\right)^{n}. $$

Collecting the sum of terms of bidegree $(m, n)$ in $(z, \overline{z})$ with $T = m + n$ in (\ref{b1}), it follows that
$$\displaystyle\sum_{m+n=T} \left(a'_{m,n}{z'}^{m}\overline{{z}'}^{n}-a_{m,n}z^{m}\overline{z}^{n}\right)= g_{T}\left (z, P(z, \overline{z})\right)-2\Re\left\{ P_{z}(z, \overline{z})f_{T-k_{0}+1} \left (z, P(z, \overline{z})\right)\right\} + \dots ,$$
where we have used the   notations
$$ g_{T}(z,w) = \displaystyle\sum_{ m+nk_{0}=T}
g_{m,n}z^{m}w^{n}\hspace{0.1 cm}\mbox{and}\hspace{0.1 cm} f_{T-k_{0}+1} (z,w) = \displaystyle\sum_{ m+nk_{0}+k_{0}-1=T}
f_{m,n}z^{m}w^{n},$$
and where the terms defined by ''$\dots$,, depend on $f_{k,l}$ with $k + k_{0}l + k_{0}-1 < T$ and as well on $g_{k,l}$ with $k + k_{0}l < T$.
 
\section{Construction of the partial normal form defined by the Fischer normalization space $\mathcal{S}_{N}$}
Considering the Fischer normalization condition  (\ref{cn}) on the  following sums of terms
$$\displaystyle\sum_{ m+n=T}
a'_{m,n}{z}^{m} \overline{z}^{n},\quad\mbox{for all $T\geq k_{0}+1$},$$
we   determine the pair of polynomials $\left(f_{T-k_{0}+1} (z,w),\hspace{0.1 cm} g_{T} (z,w)\right),$ using  the  Fischer Decompositions
 $ z^{T-k_{0}}zP_{z}(z, \overline{z}) = \alpha z^{T-k_{0}}P(z, \overline{z}) + R(z, \overline{z})$ and $\overline{z^{T-k_{0}} zP_{z}(z, \overline{z})} =\overline{ \alpha z^{T-k_{0}}P(z, \overline{z})} +\overline{ R(z, \overline{z})},$ 
where $zP_{z}(z, \overline{z}) = \alpha P(z, \overline{z}) + R(z, \overline{z})$, with $P^{\star}\left (R(z, \overline{z})\right) = 0$. 

 It is clear that the space previously introduced in (\ref{space}) is well defined, because 
each pure polynomial belongs to the kernel of the Fischer differential operator $P^{\star}$ defined in (\ref{pol}). Imposing the Fischer normalization conditions  (\ref{cn}), we determine the formal
transformation   (\ref{2800}), as long as $T \neq 0, k_{0}-1 \mod s$.  We observe that    $f_{0,n+1}$ and $f_{1,n}$ remain undetermined parameters,
for all $n\in\mathbb{N}^{\star}$. These parameters act on the higher bidegree terms helping us to impose the normalization conditions
  (\ref{extra}).

\section{The system of pseudo-weights}
  In order to apply
Huang-Yin's algorithm\cite{HY2} using the model $$w = P(z,\overline{z})+z^{s}+\overline{z}^{s},$$  we  use a different strategy defining a  system of
pseudo-weights as follows
 \begin{equation} \mbox{wt}_{\left(\gamma,\beta\right)}\left\{z^{\gamma}\overline{z}^{\beta}\right\}:=\begin{cases}\gamma+\beta ,
 & \quad\mbox{if $2\leq\gamma+\beta<k_{0}-1$ and $\gamma,\beta\neq 0$}  ;\cr
 s-1,&\quad\mbox{if $\gamma+\beta=k_{0}-1$ and $\gamma,\beta\neq 0$}
 , \end{cases}\quad\quad\quad
    \mbox{wt}\left\{z\right\}=1,\hspace{0.1 cm} \rm{wt}\left\{\overline{z}\right\}=\frac{s-1}{k_{0}-1}.
\end{equation}

For $\gamma+\beta> k_{0}$ such that $1\leq \beta\leq k_{0}-2$, we define the pseudo-weight as follows  
$$\mbox{wt}\left\{z^{\gamma}\overline{z}^{\beta} \right\}=\gamma-\left(k_{0}-1-\beta\right)+s-1.$$

Clearly, we  have  different ways how to define $\mbox{wt}\left\{z^{\gamma}\overline{z}^{\beta} \right\}$ for $\gamma+\beta>k_{0}$ with $\gamma,\beta \neq 0$. We define $$\mbox{wt}\left\{z^{\gamma}\overline{z}^{\beta} \right\}:=s-1+\rm{wt}\left\{z^{\gamma-\gamma_{1}}\overline{z}^{\beta-\beta_{1}}\right\},\quad\mbox{where $\gamma_{1}+\beta_{1}=k_{0}-1$, $\gamma_{1},\beta_{1}\neq 0$.} $$

We generally  define the pseudo-weight   
$$\mbox{wt}\left\{z^{\gamma}\overline{z}^{\beta} \right\}:=\min\mbox{wt}_{\left(\gamma,\beta\right)}\left\{z^{\gamma}\overline{z}^{\beta} \right\}.$$

It is clear by  definition that $\mbox{wt}\left\{\left(P(z,\overline{z})\right)^{2}\right\}=2s$ and $\mbox{wt}\left\{P_{z}(z,\overline{z})\right\}=s-1 $,  because $\mbox{wt}\left\{P(z,\overline{z})\right\}=s,$  but it generally holds that $\mbox{wt}\left\{z^{\gamma}\overline{z}^{\beta}\right\}\neq \rm{wt}\left\{z^{\gamma}\right\}+\rm{wt}\left\{\overline{z}^{\beta}\right\},$  and this makes our case different than the classical case (\ref{hye}), when the Bishop invariant is vanishing.

We define   the set of the
normal weights as follows
$$\mbox{wt}_{\rm{nor}} \left\{w\right\}=k_{0},\quad
\mbox{wt}_{\rm{nor}}
\left\{z\right\}=\mbox{wt}_{\rm{nor}}\left\{\overline{z}\right\}=1.$$

If $h(z,w)$ is a formal power series with no
constant term, we write
\begin{equation}\begin{split}& h(z,w)=\displaystyle\sum_{l\geq
1}h_{\rm{nor}}^{\left(l\right)}(z,w),\quad\mbox{where}\hspace{0.1 cm}\mbox{$h_{\rm{nor}}^{(l)}\left(tz,t^{k_{0}}w\right)=t^{l}h_{\rm{nor}}^{(l)}(z,w)$},\\&
\mbox{$h_{\geq l}(z,w)=\displaystyle\sum_{k\geq
l}h_{\rm{nor}}^{\left(k\right)}(z,w)$}. \end{split}\label{2.11}\end{equation}

Throughout this note, we denote by $\Theta_{N}^{\Lambda}(z,\overline{z})$ a formal power series with terms in $z,\overline{z}$ of degree at least $\Lambda$ and pseudo-weight at least $N$, and respectively by 
$\mathbb{P}_{N}^{\Lambda}(z,\overline{z})$ a bihomogeneous polynomial with terms in $z,\overline{z}$ degree   $\Lambda$ and pseudo-weight at least $N$.

\subsection{Construction Strategy}By Section $3$, we can assume that the surface $M$ is given by the  following partial normal form
\begin{equation} w=P(z,\overline{z})+\displaystyle\sum_{m+n
\geq k_{0}+1}^{T+1}a_{m,n}(z,\overline{z})+\rm{O}\left(T+2\right),
\label{ec4}\end{equation} where the Fischer normalization conditions (\ref{cn}) are satisfied, for all $k_{0}+1\leq m+n\leq T$.

We make induction on $T\geq k_{0}+1$ and we apply Huang-Yin's algorithm \cite{HY2} in order to track and compute  the parameters left undetermined  by the normalization conditions   (\ref{cn}). More exactly, when
$T+1\not\in \left\{ts;\hspace{0.1 cm}
t\in\mathbb{N}^{\star}-\left\{1,2\right\}
\right\}\cup\left\{ts+k_{0}-1;\hspace{0.1
cm}t\in\mathbb{N}^{\star}\right\},$ we impose the normalization conditions   (\ref{cn}). In the case $T+1\in\left\{ts;\hspace{0.1 cm
t\in\mathbb{N}^{\star}}-\{1\}\right\}\cup\left\{ts+k_{0}-1;\hspace{0.1
cm}t\in\mathbb{N}^{\star}\right\},$   we search a formal map which
sends our surface $M$ into a new surface $M'$ given by
\begin{equation}w'=P\left(z',\overline{z'}\right)+\displaystyle\sum_{m+n\geq k_{0}+1}^{T+1}a'_{m,n}{ z'}^{n}\overline{z'}^{m}+\rm{O}\left(T+2\right)
,\label{em2}\end{equation} where
 the Fischer normalization conditions   (\ref{cn}) and  (\ref{extra}) are satisfied, for all
$k=s+1,\dots,T$ with $k=0,k_{0}-1\hspace{0.1 cm}\rm{mod}\hspace{0.1
cm}(s)$.

 \section{Proof of Theorem \ref{teo}}

\subsection{Proof of Theorem \ref{teo}-\mbox{Case $T+1=ts+k_{0}-1$, $t\geq1$}}Throughout  this subsection, we assume $T+1=ts+k_{0}-1$, where $t\geq1$. We are looking for a biholomorphic transformation of
the following type
 \begin{equation}\left(z',w'\right)=\left(z+f(z,w),\hspace{0.1 cm}w+g(z,w)\right)=\left(z+\displaystyle\sum_{l=0}^{ts-k_{0}t}f_{\rm{nor}}^{\left(k_{0}t+l\right)}(z,w),\hspace{0.1 cm}w+\displaystyle\sum_{\tau=0}^{ts-k_{0}t}g_{\rm{nor}}^{\left(k_{0}t+k_{0}-1+\tau\right)}(z,w)\right),
\label{t1}
\end{equation} that maps $M$ into $M'$ up to  degree $T+1=ts+k_{0}-1$. 

In order for the
preceding mapping to be uniquely determined, we assume that  
$f_{\rm{nor}}^{\left(k_{0}t+l\right)}(z,w)$  has no free parameters, 
for all $l=0,\dots ts-k_{0}t$.  Substituting (\ref{t1}) into (\ref{em2}) and   simplifying using (\ref{ec4}), we obtain  
 \begin{equation}\begin{split} &\quad\quad\quad\quad\quad\quad\hspace{0.1 cm}  \displaystyle\sum_{\tau=0}^{ts-k_{0}t}g_{\rm{nor}}^{\left(k_{0}t+k_{0}-1+\tau\right)}\left(z,P(z,\overline{z})+a_{ \geq k_{0}+1}(z,\overline{z})\right)= 2\Re\left\{P_{z}(z,\overline{z})
 \displaystyle\sum_{l=0}^{ts-k_{0}t}f_{\rm{nor}}^{\left(k_{0}t+l\right)}\left(z,P(z,\overline{z})+a_{\geq k_{0}+1}(z,\overline{z})\right)\right\}\\&+2\Re\left\{\displaystyle\sum_{\alpha+\beta=2}^{k_{0}}P_{z^{\alpha}\overline{z}^{\beta}}(z,\overline{z})
 \left( \displaystyle\sum_{l=0}^{ts-k_{0}t}f_{\rm{nor}}^{\left(k_{0}t+l\right)}\left(z,P(z,\overline{z})+a_{\geq k_{0}+1}(z,\overline{z})\right)\right)^{\alpha}\left(\overline{ \displaystyle\sum_{l=0}^{ts-k_{0}t}f_{\rm{nor}}^{\left(k_{0}t+l\right)}\left(z,P(z,\overline{z})+a_{\geq k_{0}+1}(z,\overline{z})\right)}\right)^{\beta}\right\}\\&+
a'_{\geq k_{0}+1}\left(z+\displaystyle\sum_{l=0}^{ts-k_{0}t}f_{\rm{nor}}^{\left(k_{0}t+l\right)}\left(z,P(z,\overline{z})+a_{\geq k_{0}+1}(z,\overline{z})\right), \overline{z+\displaystyle\sum_{l=0}^{ts-k_{0}t}f_{\rm{nor}}^{\left(k_{0}t+l\right)}\left(z,P(z,\overline{z})+a_{\geq k_{0}+1}(z,\overline{z})\right)}\right)
-a_{\geq k_{0}+1}(z,\overline{z}).
 \end{split}\label{ecg4}
\end{equation}

Since   $f(z,w)$ does not have components of normal weight less than $k_{0}t$, collecting in (\ref{ecg4}) the sum of terms of bidegree
$(m, n)$ in $(z, \overline{z})$ with $m+n< k_{0}t+k_{0}-1$, we obtain that $a'_{
m,n}=a_{m,n}$. Collecting the sum of terms of bidegree $(m, n)$ in $(z, \overline{z})$ with $m + n = k_{0}t+ k_{0}-1$ from (\ref{ecg4}), we prove 
\bl\label{L3.1} $g^{\left(k_{0}t+k_{0}-1\right)}_{\rm{nor}} (z,w) = 0$ and $f^{\left(k_{0}t\right)}_{\rm{nor}} (z,w)=\overline{a\alpha} w^{t}-a z^{k_{0}}w^{t-1}$, where $\alpha$ is defined by (\ref{alpha}).\el
\begin{proof} Considering the corresponding iterated Fischer normalization conditions defined in (\ref{cn}), we study the terms containing $a$   $g_{k_{0}-1,t}z^{k_{0}-1} \left(P\left(z,\overline{z}\right) \right)^{t}-\left(f_{k_{0},t-k_{0}}z^{k_{0}}P_{z}\left(z,\overline{z}\right)
+ \overline{f_{0,t}P_{z}\left(z,\overline{z}\right)} P \left(z,\overline{z}\right)\right)\left(P\left(z,\overline{z}\right)\right)^{t-1}$ and $\overline{\left(f_{k_{0},t-k_{0}}z^{k_{0}}P_{z}\left(z,\overline{z}\right)
+ \overline{f_{0,t}P_{z}\left(z,\overline{z}\right)} P \left(z,\overline{z}\right)\right)}\left(P\left(z,\overline{z}\right)\right)^{t-1}$.  We obtain   $g_{k_{0}-1,t}=0$ and as well our conclusion by taking $a=-f_{k_{0},t-k_{0}}$ and by using the  
decomposition   (\ref{alpha}).  
\end{proof}

 By (\ref{2.11}) we write 
$f(z,w)=f_{\rm{nor}}^{\left(k_{0}t\right)}(z,w)+f_{\geq k_{0}t+1}(z,w)$ and
$g(z,w)=g_{\geq k_{0}t+ k_{0}}(z,w)$, where $$f_{\geq
k_{0}t+1}(z,w)=\displaystyle\sum_{k+k_{0}l\geq k_{0}t+1}f_{k,l}z^{k}w^{l}.$$

By Lemma \ref{L3.1}, it follows
that $\mathrm{wt}\left\{f_{\geq
k_{0}t+1}(z,w)\right\}\geq\displaystyle\min_{k+k_{0}l\geq k_{0}t+1}\{k+ls\}
\geq\displaystyle\min_{k+k_{0}l\geq k_{0}t+1}\{k+k_{0}l\}\geq k_{0}t+1$. In particular, 
    we obtain  
 \begin{equation}  \mathrm{wt}\left\{f_{\geq k_{0}t+1}(z,w) \right\},\mathrm{wt}\left\{\overline{f_{\geq k_{0}t+1}(z,w)} \right\}\geq k_{0}t+1
,\quad \mathrm{wt}\left\{f_{\rm{nor}}^{\left(k_{0}t\right)}(z,w)\right\},\mathrm{wt}\left\{\overline{f_{\rm{nor}}^{\left(k_{0}t\right)}(z,w)}\right\}\geq ts+k_{0}-s, \label{1e1}\end{equation}
where $w$ satisfies (\ref{ec4}).   Furthermore, by (\ref{1e1}) it follows that

\begin{equation}\begin{split} &P\left(z+\displaystyle\sum_{l=0}^{ts-k_{0}t}f_{\rm{nor}}^{\left(k_{0}t+l\right)}(z,w),\overline{z+\displaystyle\sum_{l=0}
^{ts-k_{0}t}f_{\rm{nor}}^{\left(k_{0}t+l\right)}(z,w)}\right)
=P(z,\overline{z})+2\Re\left\{P_{z}(z,\overline{z})\displaystyle\sum_{l=0}^{ts-k_{0}t}f_{\rm{nor}}^{\left(k_{0}t+l\right)}(z,w)\right\}\\&\quad\quad\quad\quad\quad\quad\quad\quad\quad\quad     +
2\Re\left\{\displaystyle\sum_{\gamma+\beta=2}^{k_{0}}P_{z^{\gamma}\overline{z}^{\beta}}(z,\overline{z})\left(
\displaystyle\sum_{l=0}^{ts-k_{0}t}f_{\rm{nor}}^{\left(k_{0}t+l\right)}(z,w)\right)^{\gamma}\left(\overline{
\displaystyle\sum_{l=0}^{ts-k_{0}t}f_{\rm{nor}}^{\left(k_{0}t+l\right)}(z,w)}\right)^{\beta}\right\}\\&\quad \quad \quad \quad \quad \quad \quad 
=P(z,\overline{z})+2\Re\left\{P_{z}(z,\overline{z})f_{\rm{nor}}^{\left(k_{0}t\right)}(z,w)\right\}+2\Re\left\{\left(P_{z}(z,\overline{z})+\Theta_{s}^{k_{0}}(z,\overline{z})\right)f_{\geq k_{0}t+1}(z,w)\right\}
+\Theta_{ts+k_{0}}^{k_{0}t+k_{0}}(z,\overline{z}),\end{split}\label{ecg444}\end{equation}
where $w$ satisfies (\ref{ec4}) and $\rm{wt}\left\{\overline{\Theta_{ts+k_{0}}^{k_{0}t+k_{0}}(z,\overline{z})}\right\}\geq ts+k_{0}$.

 In order to track the action of the free parameter $a$, we   prove
 \bl\label{L3.3} For all  $m,n\geq 1$ with $m+n\geq k_{0}+1$ and $w$ satisfying (\ref{ec4}), we
have 
\begin{equation} \left(z+f(z,w)\right)^{m}\left(\overline{z+f(z,w)}\right)^{n}=z^{m}\overline{z}^{n}
+2\Re\left\{\Theta_{s}^{k_{0}}(z,\overline{z})
 f_{\geq
k_{0}t+1}(z,w)\right\}+\Theta_{ts+k_{0}}^{k_{0}t+k_{0}}(z,\overline{z}),\end{equation}
where
$\mathrm{wt}\left\{\overline{\Theta_{ts+k_{0}}^{k_{0}t+k_{0}}(z,\overline{z})}\right\}\geq
 ts+k_{0}$.
 \el
\begin{proof} By the Taylor expansion, we obtain  
\begin{equation}\begin{split}\left(z+f(z,w)\right)^{m}\left(\overline{z+f(z,w)}\right)^{n}&=z^{m}\overline{z}^{n}
+2\Re\left\{\Theta_{s}^{k_{0}}(z,\overline{z})\left(f_{\rm{nor}}^{\left(k_{0}t\right)}(z,w)+
 f_{\geq
k_{0}t+1}(z,w)\right)\right\} \\&=z^{m}\overline{z}^{n}
+2\Re\left\{\Theta_{s}^{k_{0}}(z,\overline{z})
 f_{\geq
k_{0}t+1}(z,w)\right\}+\Theta_{ts+k_{0}}^{k_{0}t+k_{0}}(z,\overline{z}),\end{split}\end{equation}
where
$\mathrm{wt}\left\{\overline{\Theta_{ts+k_{0}}^{k_{0}t+k_{0}}(z,\overline{z})}\right\}\geq
 ts+k_{0}$. \end{proof}
 
Analogously, we obtain  
 \bl\label{L3.4} For $w$ satisfying (\ref{ecg4}) and for all $k>s$, we have  
 \begin{equation} \left(z+f(z,w)\right)^{k}=z^{k}+2\Re \left\{
\Theta_{s}^{k_{0}}(z,\overline{z})f_{\geq
k_{0}t+1}(z,w)\right\}+\Theta_{ts+k_{0}}^{k_{0}t+k_{0}}(z,\overline{z}),\end{equation}
where
$\mathrm{wt}\left\{\overline{\Theta_{ts+k_{0}}^{k_{0}t+k_{0}}(z,\overline{z})}\right\}\geq
 ts+k_{0}$.
 \el
  \bl\label{L3.5} For $f_{\rm{nor}}^{\left(k_{0}t\right)}(z,w)$
given by Lemma \ref{L3.1}  and $w$ satisfying (\ref{ecg4}) we have
\begin{equation}2\Re\left\{P\left ( f_{\rm{nor}}^{\left(k_{0}t\right)}\left(z,w\right),z\right
)\right\}=2\Re\left\{ a\alpha z^{s-1+k_{0}}w^{t-1}\right\}-2\Re\left\{az^{k_{0}-1}R(z,\overline{z})w^{t-1}\right\} +\Theta_{ts+k_{0}}^{k_{0}t+k_{0}}(z,\overline{z}),\label{90q}\end{equation}
where
$\mathrm{wt}\left\{\overline{\Theta_{ts+k_{0}}^{k_{0}t+k_{0}}(z,\overline{z})}\right\}\geq
 ts+k_{0}$.
 \el
\begin{proof} Because $w^{t}=\overline{w}^{t}+w^{t}-\overline{w}^{t}=\overline{w}^{t}+\Theta_{ts-s+s+1}^{tk_{0}+1}(z,\overline{z})$, while
\begin{equation}2\Re\left\{P\left ( f_{\rm{nor}}^{\left(k_{0}t\right)}\left(z,w\right),z\right)
\right\}=2\Re\left\{\left(\overline{a\alpha}w-az^{k_{0}-1}z\right)P_{z} (z,\overline{z})w^{t-1}\right\}
+\Theta_{ts+k_{0}}^{k_{0}t+k_{0}}(z,\overline{z}),\end{equation}
we get  (\ref{90q}) using  the  decomposition  (\ref{alpha}). Because 
$\partial_{z}\left(P(z,\overline{z})\right)=\overline{z}^{k_{0}-1}+Q(z,\overline{z})$, we have 
 $\mbox{wt}\left\{Q(z,\overline{z}) w^{t}\right\},\mbox{wt}\left\{\overline{Q(z,\overline{z}) w^{t}}\right\} \geq ts+k_{0}$. 
 \end{proof}

 \bl\label{L3.6} For $w$ satisfying  (\ref{ecg4})
we have 
\begin{equation} 2\Re\{ \left(z+f(z,w)\right)^{s}\}=2\Re\left\{z^{s}-sa z^{s-1+k_{0}}w^{t-1}\right\}+2\Re\left\{
\left(sz^{s-1}+\Theta_{s}^{k_{0}}(z,\overline{z}) \right)f_{\geq
k_{0}t+1}(z,w)\right\}+\Theta_{ts+k_{0}}^{k_{0}t+k_{0}}(z,\overline{z}), \label{24}\end{equation}
where
$\mathrm{wt}\left\{\overline{\Theta_{ts+k_{0}}^{k_{0}t+k_{0}}(z,\overline{z})}\right\}\geq
 ts+k_{0}$.\el
 \begin{proof}
  Using the Taylor expansion, it follows that
  \begin{equation*}\begin{split}\left(z+f(z,w)\right)^{s}&=2\Re\left\{z^{s}\right\}
  +2\Re\left\{\left(sz^{s-1}+\Theta_{s}^{k_{0}}(z,\overline{z}) \right)\left(f_{\rm{nor}}^{\left(k_{0}t\right)}(z,w)+f_{\geq k_{0}t+1} (z,w)\right)\right\}  \\&=2\Re\left\{z^{s}-saz^{s-1+k_{0}}w^{t-1} \right\}
  +2\Re\left\{\left(sz^{s-1}+\Theta_{s}^{k_{0}}(z,\overline{z}) \right) f_{\geq k_{0}t+1} (z,w) \right\}+\Theta_{ts+k_{0}}^{k_{0}t+k_{0}}(z,\overline{z}),\end{split}\end{equation*}
where $w$ is given by (\ref{ecg4}) and $\mathrm{wt}\left\{\overline{\Theta_{ts+k_{0}}^{k_{0}t+k_{0}}(z,\overline{z})}\right\}\geq
 ts+k_{0}$.
\end{proof}

  By Lemmas \ref{L3.1}-\ref{L3.6} and by (\ref{ecg4}), (\ref{ecg444}), we obtain  
\begin{equation}\begin{split}
g_{\geq k_{0}t+k_{0}}\left(z,w)\right)=&2 \Re\left\{\left(\left( \alpha -s\right) a z^{s+k_{0}-1}\right)w^{t-1}\right\}
 +2\Re\left\{\left( P_{z}(z,\overline{z})+sz^{s-1}+\Theta_{s}^{k_{0}}(z,\overline{z})\right)f_{\geq
k_{0}t+1} (z,w)\right\}  \\&+a'_{\geq k_{0}t+k_{0}}(z,\overline{z})-a_{\geq k_{0}t+k_{0}}(z,\overline{z})+\Theta_{ts+k_{0}}^{k_{0}t+k_{0}}(z,\overline{z})-2\Re\left\{ az^{k_{0}-1}R(z,\overline{z})\left(w^{t-1}-\left(P(z,\overline{z})\right)^{t-1}\right)\right\} ,
\end{split}  \label{27}\end{equation}
where
$\mathrm{wt}\left\{\overline{\Theta_{ts+k_{0}}^{k_{0}t+k_{0}}(z,\overline{z})}\right\}\geq
 ts+k_{0}$, $w$ satisfies (\ref{ec4}) and $\alpha$ is defined by (\ref{alpha}).

 Assume  $t=1$. Collecting the terms of total degree $k\leq s+k_{0}-1$ in $(z,\overline{z})$  in (\ref{27}) we
 find the  pair of polynomials
$$\left(f_{\rm{nor}}^{\left(k-k_{0}+1\right)}(z,w),\hspace{0.1 cm}g_{\rm{nor}}^{\left(k\right)}(z,w)\right).$$

 Collecting the terms of total degree $m+n=s+k_{0}-1$ in
$(z,\overline{z})$ in (\ref{27}), we  obtain
\begin{equation}\begin{split}  g_{\rm{nor}}^{\left(s+k_{0}-1\right)}\left(z,P(z,\overline{z})\right)=&2\Re\left\{\left( \alpha -s\right) a  z^{s+k_{0}-1}  \right\} +2\Re\left\{P_{z}(z,\overline{z})f_{\rm{nor}}^{(s)}\left(z,P(z,\overline{z})\right)\right\}
  \\&+a'_{s+k_{0}-1}(z,\overline{z})-a_{s+k_{0}-1}(z,\overline{z})+\left(\Theta_{1}\right)_{s+k_{0}}^{k_{0}+k_{0}}(z,\overline{z}).\end{split}
  .\label{28}
\end{equation}

Imposing the normalization condition $a'_{0,s+k_{0}-1}=0$, we compute the parameter $a$.

 Assuming  $t\geq 2$, we prove:
\bl \label{L3.7}Let $N_{s}:=ts+k_{0}$. For all $0\leq j\leq t-1$ and
$p\in\left[k_{0}t+j\left(s-k_{0}\right)+k_{0},k_{0}t+(j+1)\left(s-k_{0}\right)+k_{0}-1\right]$, we have 
\begin{equation}\begin{split} g_{\geq p}(z,w)=& 2\Re\left\{z^{s(j+1)+k_{0}-1}
\left( \left( \alpha -s\right) a\left(1-\frac{s}{\alpha}\right)^{j}w^{t-j-1}-R(z,\overline{z})\sum_{\beta+\gamma=t-j-2}v_{\beta\gamma}w^{\beta}\left(P(z,\overline{z})\right)^{\gamma}\right) \right\}\\&+2\Re\left\{\left( P_{z}(z,\overline{z})+
sz^{s-1}+\Theta_{s}^{k_{0}}(z,\overline{z})\right)f_{\geq
p-k_{0}+1}(z,w) \right\}   +a'_{\geq
p}(z,\overline{z})-a_{\geq
p}(z,\overline{z})+\Theta_{N_{s}}^{p}(z,\overline{z}) ,\end{split}
\label{30}
\end{equation}
where
$\mathrm{wt}\left\{\overline{\Theta_{N_{s}}^{p}(z,\overline{z})}\right\}\geq
 N_{s}$, $w$ satisfies (\ref{ec4}) and $\alpha$ is defined by (\ref{alpha}).
  \el
\begin{proof} We organize our proof in two steps.

$\bf{Step\hspace{0.1 cm}1.}$ When $s=k_{0}+1$ this step is obvious.
Assume that $s>k_{0}+1$. Let  $p_{0}=k_{0}t+j\left(s-k_{0}\right)+k_{0}$, where
$j\in\left[0,t-1\right]$.
 We make induction on $p\in\left[k_{0}t+j\left(s-k_{0}\right)+k_{0}, k_{0}t+(j+1)\left(s-k_{0}\right)+k_{0}-1\right]$. For $j=0$
  the lemma is satisfied (see  
(\ref{27})). Let $p\geq p_{0}$ such that $p+1\leq
k_{0}t+(j+1)\left(s-k_{0}\right)+k_{0}-1$. Collecting the terms of bidegree $(m,n)$ in
$(z,\overline{z})$ in (\ref{30}) with $m+n=p$, we obtain  
\begin{equation}\hspace{0.1 cm}g_{\rm{nor}}^{(p)}\left(z,P(z,\overline{z})\right)=
2\Re\left\{ P_{z}(z,\overline{z})f_{\rm{nor}}^{\left(p-k_{0}+1\right)}\left(z,P(z,\overline{z})\right)\right\}+
a'_{p}(z,\overline{z})-a_{p}(z,\overline{z})+\mathbb{P}_{N_{s}}^{p}(z,\overline{z}).\label{A}\end{equation}

Because the Fischer normalization conditions   (\ref{cn}) are satisfied,  we  find a normalized solution
 $\left(f_{\rm{nor}}^{\left(p-k_{0}+1\right)}(z,w),g_{\rm{nor}}^{(p)}(z,w)\right)$  for
(\ref{A}) and   as in the case of Huang-Yin\cite{HY2}, the following estimates hold
\begin{equation}\begin{split}&  \left\{g_{\rm{nor}}^{(p)}(z,w)\right\} ,\hspace{0.1
cm}\mathrm{wt}\left\{g_{\rm{nor}}^{(p)}(z,w)-g_{\rm{nor}}^{(p)}\left(z,P(z,\overline{z})\right)\right\}\geq N_{s}, \quad \mathrm{wt}\left\{f_{\rm{nor}}^{\left(p-k_{0}+1\right)}(z,w)\right\},\hspace{0.1
cm}\mathrm{wt}\left\{\overline{f_{\rm{nor}}^{\left(p-k_{0}+1\right)}\left(z,w\right)}\right\}, \\& \mathrm{wt}\left\{\overline{f_{\rm{nor}}^{\left(p-k_{0}+1\right)}(z,w)-f_{\rm{nor}}^{\left(p-k_{0}+1\right)}\left(z,P(z,\overline{z})\right)}\right\} ,
\quad\mathrm{wt}\left\{f_{\rm{nor}}^{\left(p-k_{0}+1\right)}(z,w)-f_{\rm{nor}}^{\left(p-k_{0}+1\right)}\left(z,P(z,\overline{z})\right)\right\}\geq N_{s}-s+1,\end{split}
  \label{31345}\end{equation}
where $w$ is given by (\ref{ec4}). By (\ref{31345}) we obtain  
\begin{equation}\begin{split}&\quad\quad\quad\quad\quad\quad\quad\quad\quad\quad\quad\quad\quad\quad \quad\quad\quad  g_{\rm{nor}}^{(p)}(z,w)-g_{\rm{nor}}^{(p)}\left(z,P(z,\overline{z})\right)=\Theta_{N_{s}}^{p+1}(z,\overline{z})',\\&
 2\Re\left\{P_{z}(z,\overline{z})\left(f_{\rm{nor}}^{\left(p-k_{0}+1\right)}(z,w)-f_{\rm{nor}}^{\left(p-k_{0}+1\right)}\left(z,P(z,\overline{z})\right)\right)
\right\}=\Theta_{N_{s}}^{p+1}(z,\overline{z})',\quad\quad 2\Re\left\{\Theta_{s}^{k_{0}}(z,\overline{z}) f_{\rm{nor}}^{\left(p-k_{0}+1\right)}(z,w)\right\}=\Theta_{N_{s}}^{p+1}(z,\overline{z})', \end{split} \label{400} \end{equation}
and each of the preceding formal power series
$\Theta_{N_{s}}^{p+1}(z,\overline{z})'$ has the property
$\mathrm{wt}\left\{\overline{\Theta_{N_{s}}^{p+1}(z,\overline{z})'}\right\}\geq
 N_{s}$.  
 
  Substituting  $f_{\geq  p-k_{0}+1 }(z,w)=f_{\rm{nor}}^{\left(p-k_{0}+1\right)}(z,w)+f_{\geq
 p-k_{0}+2 }(z,w)$ and $g_{\geq p}(z,w)=g_{\rm{nor}}^{(p)}(z,w)+g_{\geq
p+1}(z,w)$ into (\ref{30}), we obtain
 \begin{equation}\begin{split}  g_{\rm{nor}}^{\left(p\right)}(z,w)+g_{\geq p+1}(z,w)=& 2\Re\left\{z^{s(j+1)+k_{0}-1}
\left( \left( \alpha -s\right) a\left(1-\frac{s}{\alpha}\right)^{j}w^{t-j-1}- R(z,\overline{z})\sum_{\beta+\gamma=t-j-2}v_{\theta\gamma}w^{\beta}\left(P(z,\overline{z})\right)^{\gamma}\right) \right\}  
\\&+2\Re\left\{\left( P_{z}(z,\overline{z})+
sz^{s-1}+\Theta_{s}^{k_{0}}(z,\overline{z})\right)\left(f_{\rm{nor}}^{\left(p-k_{0}+1\right)}(z,w)+f_{\geq
 p-k_{0}+2 }(z,w)\right)\right\}\\&+a'_{\geq p}(z,\overline{z})-a_{\geq
p}(z,\overline{z})+\mathbb{P}_{N_{s}}^{p}(z,\overline{z})
 +\Theta_{N_{s}}^{p+1}(z,\overline{z}) .\end{split}  \label{32}\end{equation}

 By making a simplification in (\ref{32}) by using
(\ref{A}),  it follows that
\begin{equation}\begin{split} g_{\geq
p+1}(z,w)=&2\Re\left\{z^{s(j+1)+k_{0}-1}
\left( \left( \alpha -s\right) a\left(1-\frac{s}{\alpha}\right)^{j}w^{t-j-1}- R(z,\overline{z})\sum_{\beta+\gamma=t-j-2}v_{\theta\gamma}w^{\beta}\left(P(z,\overline{z})\right)^{\gamma}\right) \right\}\\&+2\Re\left\{ \left(P_{z}(z,\overline{z} )+
 sz^{s-1}+\Theta_{s}^{k_{0}}(z,\overline{z})\right) f_{\geq
p-k_{0}+2}(z,w)\right\} +a'_{\geq
p+1}(z,\overline{z})-a_{\geq
p+1}(z,\overline{z})+J(z,\overline{z}) +\Theta_{N_{s}}^{p+1}(z,\overline{z}) ,\end{split} \label{33}
\end{equation}
where
 we have used the following notation
 \begin{equation}\begin{split}J(z,\overline{z})=&2\Re \left\{P_{z}(z,\overline{z})\left(f_{\rm{nor}}^{\left(p-k_{0}+1\right)}(z,w)-f_{\rm{nor}}^{\left(p-k_{0}+1\right)}\left(z,P(z,\overline{z})\right)\right)\right\}+2\Re \left\{\left(
sz^{s-1}+\Theta_{s}^{k_{0}}(z,\overline{z})\right)
 f_{\rm{nor}}^{\left(p-k_{0}+1\right)}(z,w)\right\}\\&+g_{\rm{nor}}^{(p)}\left(z,P(z,\overline{z})\right)-g_{\rm{nor}}^{(p)}(z,w).
\end{split}  \label{VB}\end{equation}
  
By   (\ref{31345}) and (\ref{400}) it follows that
$J(z,\overline{z})=\Theta_{N_{s}}^{p+1}(z,\overline{z})$, where
$\mathrm{wt}\left\{\overline{\Theta_{N_{s}}^{p+1}(z,\overline{z})}\right\}\geq
 N_{s}$.

 $\bf{Step\hspace{0.1 cm}2.}$ Assume that
we have proved Lemma \ref{L3.7} for
$p\in\left[k_{0}t+j\left(s-k_{0}\right)+k_{0},k_{0}t+(j+1)\left(s-k_{0}\right)+k_{0}-1\right]$, for 
$j\in[0,t-1]$. We   prove  Lemma \ref{L3.7} for
$p\in\left[k_{0}t+(j+1)\left(s-k_{0}\right)+k_{0},k_{0}t+(j+2)\left(s-k_{0}\right)+k_{0}-1\right]$. Adapting Huang-Yin's strategy\cite{HY2}, we define  
 \begin{equation}\begin{split}& f_{\rm{nor}}^{\left(\Lambda \right)}(z,w)=f_{1}^{\left(\Lambda \right)}(z,w)+f_{2}^{\left(\Lambda \right)}(z,w),\\& f_{1}^{\left(\Lambda\right)}(z,w)=-\frac{ a }{\alpha}\left( \alpha -s\right)\left(1-\frac{s}{\alpha}\right)^{j} z^{s(j+1)+k_{0} }
w^{t-j-2},\quad\mbox{for $\Lambda= k_{0}t+(j+1)\left(s-k_{0}\right) $.}\end{split} \label{35}\end{equation}

Substituting
$f_{\geq\Lambda}(z,w)=f_{\rm{nor}}^{\left(\Lambda \right)}(z,w)+f_{\geq\Lambda+1}(z,w)$
and
$g_{\geq\Lambda+k_{0}-1}(z,w)=g_{\rm{nor}}^{\left(\Lambda+k_{0}-1\right)}(z,w)+g_{\Lambda+k_{0}}(z,w)$
in (\ref{30}), we obtain  

\begin{equation}\begin{split}  g_{\rm{nor}}^{\left(\Lambda+k_{0}-1\right)}(z,w)+g_{\geq\Lambda+k_{0}}(z,w)=&2\Re\left\{z^{s(j+1)+k_{0}-1}
\left( \left( \alpha -s\right) a\left(1-\frac{s}{\alpha}\right)^{j}w^{t-j-1}-R(z,\overline{z})\sum_{\beta+\gamma=t-j-2}v_{\theta\gamma}w^{\beta}\left(P(z,\overline{z})\right)^{\gamma}\right) \right\} \\&+2\Re\left\{
\left(P_{z}(z,\overline{z})+ sz^{s-1}
+\Theta_{s}^{k_{0}}(z,\overline{z})\right)\left(f_{\rm{nor}}^{\left(\Lambda \right)}(z,w)+f_{\geq
\Lambda+1}(z,w)\right)\right\}
 \\& 
+a'_{\geq\Lambda+k_{0}-1}(z,\overline{z})-a_{\geq\Lambda+k_{0}-1}(z,\overline{z})+\mathbb{P}_{N_{s}}^{\Lambda+k_{0}-1}(z,\overline{z})
+\Theta_{N_{s}}^{\Lambda+k_{0}}(z,\overline{z}) .\end{split}\label{38}\end{equation}

Rewriting  (\ref{38}) using (\ref{35}),  we obtain  
\begin{equation}\begin{split} g_{\rm{nor}}^{\left(\Lambda+k_{0}-1\right)}(z,w)+g_{\geq\Lambda+k_{0}}(z,w)=&
2\Re\left\{\left(P_{z}(z, \overline{z})+
sz^{s-1}+\Theta_{s}^{k_{0}}(z,\overline{z})\right)\left(f_{\geq\Lambda+1}(z,w)+f^{\left(\Lambda\right)}_{2}(z,w)\right)\right\}
\\&  +a'_{\geq\Lambda+k_{0}-1}(z,\overline{z})-a_{\geq\Lambda+k_{0}-1}(z,\overline{z})
+\Theta_{N_{s}}^{\Lambda+k_{0}}(z,\overline{z})+J(z,\overline{z}),\end{split} \label{39v}\end{equation}
where we have used the following notation
\begin{equation*}\begin{split} J(z,\overline{z})&=2\Re\left \{P_{z}(z,\overline{z})\left(f^{\left(\Lambda \right)}_{1}(z,w)-f^{\left(\Lambda \right)}_{1}\left(z,P(z,\overline{z})\right) \right)\right\}+2\Re \left\{\left(\Theta_{s}^{k_{0}}(z,\overline{z})+sz^{s-1}\right)f_{1}^{\left(\Lambda\right)}(z,w)\right\}\\& +2\Re\left\{z^{s(j+1)+k_{0}-1}
\left( \left( \alpha -s\right) a\left(1-\frac{s}{\alpha}\right)^{j}\left(w^{t-j-1}-\left(P(z,\overline{z})\right)^{t-j-1}\right)\right.\right. \\& \left.\left. \quad\quad \quad\hspace{0.21 cm} -R(z,\overline{z})\sum_{\beta+\gamma=t-j-2}v_{\beta\gamma}\left(w^{\beta}-\left(P(z,\overline{z})\right)^{\beta}\right)\left(P(z,\overline{z})\right)^{\gamma}\right) \right\}. \end{split}   \end{equation*}

Because the Fischer normalization conditions   (\ref{cn}) remain preserved, it follows  that
$$J(z,\overline{z})= 2\Re\left\{z^{s(j+2)+k_{0}-1}
\left( \left( \alpha -s\right) a\left(1-\frac{s}{\alpha}\right)^{j+1}w^{t-j-2}- R(z,\overline{z})\sum_{\beta+\gamma=t-j-3}v_{\beta\gamma}w^{\beta}\left(P(z,\overline{z})\right)^{\gamma}\right) \right\}+\Theta_{N_{s}}^{\Lambda+k_{0}}(z,\overline{z}) .$$

Collecting the
terms of bidegree $(m,n)$ in $(z,\overline{z})$ in (\ref{39v}) with
$m+n=\Lambda+k_{0}-1$,  we obtain  
\begin{equation}  g_{\rm{nor}}^{\left(\Lambda+k_{0}-1\right)}\left(z,P(z,\overline{z})\right)=  2\Re\left\{
P_{z}(z,\overline{z})f_{2}^{\left(\Lambda\right)}\left(z,P(z,\overline{z}\right)\right\} +a'_{\Lambda+k_{0}-1}(z,\overline{z})-a_{\Lambda+k_{0}-1}(z,\overline{z})
+\mathbb{P}_{N_{s}}^{\Lambda+k_{0}-1}(z,\overline{z}).
 \label{34xxx} \end{equation}
 
 Because the corresponding Fischer normalization conditions   (\ref{cn}) are satisfied, we find a solution
 $\left(g_{\rm{nor}}^{\left(\Lambda+k_{0}-1\right)}(z,w),\hspace{0.1 cm}f_{2}^{\left(\Lambda \right)}(z,w)\right)$ 
for (\ref{34xxx}), satisfying the following estimates
\begin{equation}\begin{split}&\mathrm{wt}\left\{g_{\rm{nor}}^{\left(\Lambda+k_{0}-1\right)}(z,w)-g_{\rm{nor}}^{\left(\Lambda+k_{0}-1\right)}\left(z,P(z,\overline{z})\right)\right\},
\mathrm{wt}\left\{g_{\rm{nor}}^{\left(\Lambda +k_{0}-1\right)}(z,w)\right\},\hspace{0.1
cm}\mathrm{wt}\left\{ g_{\rm{nor}}^{\left(\Lambda+k_{0}-1\right)}\left(z,P(z,\overline{z})\right) \right\} \geq N_{s},\\& \quad\quad\quad\quad\quad\quad\quad\quad\hspace{0.2 cm}\mathrm{wt}\left\{\overline{f_{2}^{\left(\Lambda\right)}(z,w})\right\}, \mathrm{wt}\left\{f_{2}^{\left(\Lambda \right)}(z,w)\right\} ,
\mathrm{wt}\left\{f_{2 }^{\left(\Lambda \right)}\left(z,P(z,\overline{z})\right)\right\},\mathrm{wt}\left\{\overline{f_{2 }^{\left(\Lambda \right)}\left(z,P(z,\overline{z})\right)}\right\}\geq N_{s}-s+1, \end{split} \label{31}\end{equation}
 where  $w$ satisfies
(\ref{ec4}). As a consequence of (\ref{31}), we obtain  
\begin{equation}\begin{split}& 2\Re\left\{\left(
f_{2}^{\left(\Lambda \right)}(z,w)-f_{2}^{\left(\Lambda \right)}\left(z,P(z,\overline{z})\right)\right)P_{z}(z,\overline{z})\right\}=\Theta_{N_{s}}^{\Lambda+k_{0}}(z,\overline{z}),  \quad\quad  
2\Re\left\{\left( sz^{s-1}+\Theta_{s}^{k_{0}}(z,\overline{z})\right)f_{2}^{\left(\Lambda \right)}(z,w)\right\}=\Theta_{N_{s}}^{\Lambda+k_{0}}(z,\overline{z}),  \\&\quad \quad \quad \quad \quad \quad\quad \quad \quad \quad \quad  \hspace{0.1 cm} 
g_{\rm{nor}}^{\left(\Lambda+k_{0}-1\right)}(z,w)-g_{\rm{nor}}^{\left(\Lambda+k_{0}-1\right)}\left(z,P(z,\overline{z})\right)=\Theta_{N_{s}}^{\Lambda+k_{0}}(z,\overline{z}),
\end{split}
 \label{37}
\end{equation}
where $w$ satisfies (\ref{ec4}) and each of the preceding formal
power series has the property
$\mathrm{wt}\left\{\overline{\Theta_{N_{s}}^{\Lambda+k_{0}}(z,\overline{z})}\right\}\geq
 N_{s}$.
 
  Following the computations the first step of the proof of this lemma, we  are able to finish the proof, because the multiple of $R(z,\overline{z})$ does not appear when $t=j-1$ in  (\ref{30}).
 \end{proof}
Collecting the  terms of bidegree $(m,n)$ in $(z,\overline{z})$
with $m+n=ts+k_{0}-1$ and $t=j-1$ in  (\ref{30}), we obtain
\begin{equation}\begin{split}g_{\rm{nor}}^{\left(ts+k_{0}-1\right)}\left(z,P(z,\overline{z})\right)=&2 \Re\left\{\left( \alpha -s\right) a
\left(1-\frac{s}{\alpha}\right)^{t-1}z^{ts+k_{0}-1} \right\}+2\Re\left\{P_{z}(z,\overline{z})
f_{\rm{nor}}^{(ts)}\left(z,P(z,\overline{z})\right)\right\} \\&
+a'_{ts+k_{0}-1 }(z,\overline{z})-a_{ts+k_{0}-1 }(z,\overline{z})+\left(\Theta_{1}\right)_{N_{s}}^{ts+k_{0}-1}(z,\overline{z}),\end{split} \label{46q}\end{equation}
where $\rm{wt}\left\{\overline{\left(\Theta_{1}\right)_{N_{s}}^{ts+k_{0}-1}(z,\overline{z})}\right\}\geq N_{s}$. 

By imposing the corresponding Fischer normalization conditions  (\ref{cn}) we find the solution
$\left(f_{\rm{nor}}^{\left(ts\right)}(z,w),\hspace{0.1 cm}g_{\rm{nor}}^{\left(ts+k_{0}-1\right)}(z,w)\right)$ for
(\ref{46q}).  The parameter
$a$ is computed  by considering the corresponding Fischer normalization condition from
(\ref{extra}) in this case.   By composing the map that sends $M$ into
(\ref{ec4}) with the map (\ref{t1}) we obtain our formal
transformation that sends $M$ into $M'$ up  to degree $ts+k_{0}-1$.

\subsection{Proof of Theorem \ref{teo}-\mbox{Case $T+1=ts+s$, $t\geq1$}}Throughout this subsection, we assume
$T+1=ts+s$, where $t\geq 1$. We are looking for a biholomorphic transformation of
the following type
 \begin{equation} \left(z',w'\right)=\left(z+f(z,w),\hspace{0.1 cm}w+g(z,w)\right)=\left(z+\displaystyle\sum_{l=0}^{\left(s-k_{0}\right)\left(t+1\right)}f_{\rm{nor}}^{\left(k_{0}t+l+1\right)}(z,w),\hspace{0.1 cm}w+\displaystyle\sum_{\tau=0}^{\left(s-k_{0}\right)\left(t+1\right)}g_{\rm{nor}}^{\left(k_{0}t+k_{0}+\tau\right)}(z,w)\right),
\label{t1se}\end{equation} that maps $M$ into $M'$ up to   degree $T+1=ts+s$.

 In order for the
preceding mapping to be uniquely determined, we assume that
$f_{\rm{nor}}^{\left(k_{0}t+l+1\right)}(z,w)$  has no free parameters,
for all $l=0,\dots, ts+s-k_{0}t-k_{0}+1=\left(s-k_{0}\right)\left(t+1\right)$. Substituting (\ref{t1se}) into (\ref{em2}) and simplifying using (\ref{ec4}), we obtain  
 \begin{equation}\begin{split}  &\displaystyle\sum_{\tau=0}^{\left(s-k_{0}\right)\left(t+1\right)}g_{\rm{nor}}^{\left(k_{0}t+k_{0}+\tau\right)}\left(z,P(z,\overline{z})+a_{ \geq k_{0}+1}(z,\overline{z})\right)=  2\Re\left\{P_{z}(z,\overline{z})
 \displaystyle\sum_{l=0}^{\left(s-k_{0}\right)\left(t+1\right)}f_{\rm{nor}}^{\left(k_{0}t+l+1\right)}\left(z,P(z,\overline{z})+a_{\geq k_{0}+1}(z,\overline{z})\right)\right\}\\&\quad\quad\quad\quad\quad\quad\quad\quad\quad\quad\quad\quad\quad\quad\quad\quad+2\Re\left\{\displaystyle\sum_{\alpha+\beta=2}^{k_{0}}P_{z^{\alpha}\overline{z}^{\beta}}(z,\overline{z})
 \left( \displaystyle\sum_{l=0}^{\left(s-k_{0}\right)\left(t+1\right)}f_{\rm{nor}}^{\left(k_{0}t+l+1\right)}\left(z,P(z,\overline{z})+a_{\geq k_{0}+1}(z,\overline{z})\right)\right)^{\alpha}\right. \\& \left.\left(\overline{ \displaystyle\sum_{l=0}^{\left(s-k_{0}\right)\left(t+1\right)}f_{\rm{nor}}^{\left(k_{0}t+l+1\right)}\left(z,P(z,\overline{z})+a_{\geq k_{0}+1}(z,\overline{z})\right)}\right)^{\beta}\right\}+
a'_{\geq k_{0}+1}\left(z+\displaystyle\sum_{l=0}^{\left(s-k_{0}\right)\left(t+1\right)}f_{\rm{nor}}^{\left(k_{0}t+l+1\right)}\left(z,P(z,\overline{z})+a_{\geq k_{0}+1}(z,\overline{z})\right),\right. \\& \quad\quad\quad\quad\quad\quad\quad\quad\quad\quad\quad\quad\quad\quad\quad\quad\quad
\quad\quad\quad\quad\quad\quad  \left.\overline{z+\displaystyle\sum_{l=0}^{\left(s-k_{0}\right)\left(t+1\right)}f_{\rm{nor}}^{\left(k_{0}t+l+1\right)}\left(z,P(z,\overline{z})+a_{\geq k_{0}+1}(z,\overline{z})\right)}\right) -a_{\geq k_{0}+1}(z,\overline{z}).\end{split}
 \label{ecg4se}
\end{equation}

Since $f(z,w)$ and $g(z,w)$ do  not have components of normal weight less than $k_{0}t+1$, collecting in (\ref{ecg4}) the sum of terms of bidegree
$(m, n)$ in $(z, \overline{z})$ with $m+n< k_{0}t+k_{0}$, we obtain that $a'_{
m,n}=a_{m,n}$. Collecting the sum of terms of bidegree $(m, n)$ in $(z, \overline{z})$ with $m+n=k_{0}t+k_{0}$ from (\ref{ecg4}), we prove

\bl\label{L3.1se} $g^{\left(k_{0}t+k_{0}\right)}_{\rm{nor}} (z,w) = \left(\alpha a+\overline{\alpha a}\right)w^{t+1}$ and $f^{\left(k_{0}t+1\right)}_{\rm{nor}} (z,w)=a z  w^{t}$, where $\alpha$ is defined by (\ref{alpha}).\el
\begin{proof} Considering the corresponding iterated Fischer normalization conditions defined in (\ref{cn}), the terms that provides us the
undetermined parameter $a$ are the terms $g_{0,t+1}\left (P(z, \overline{z})\right)^{t+1}-\left(f_{1,t}zP_{z}(z, \overline{z}) + \overline{f_{1,t}zP_{z}(z, \overline{z})}\right)\left(P(z, \overline{z})\right)^{t}.$  By the uniqueness of the Fischer decomposition, we obtain   our conclusion by taking $a=f_{1,t}$ and by using the decomposition (\ref{alpha}).
\end{proof}

 By (\ref{2.11}) we write  $f(z,w)=f_{\rm{nor}}^{\left(k_{0}t+1\right)}(z,w)+f_{\geq k_{0}t+2}(z,w)$ and
$g(z,w)=g_{\geq k_{0}t+ k_{0}}(z,w)$, where $$f_{\geq
k_{0}t+2}(z,w)=\displaystyle\sum_{k+k_{0}l\geq k_{0}t+2}f_{k,l}z^{k}w^{l}.$$

By Lemma \ref{L3.1se}  it follows
that $\mathrm{wt}\left\{f_{\geq
k_{0}t+2}(z,w)\right\}\geq\displaystyle\min_{k+k_{0}l\geq k_{0}t+2}\{k+ls\}
\geq\displaystyle\min_{k+k_{0}l\geq k_{0}t+2}\{k+k_{0}l\}\geq k_{0}t+2$. In particular, 
    we obtain  
 \begin{equation}  \mathrm{wt}\left\{f_{\geq k_{0}t+2}(z,w) \right\},\mathrm{wt}\left\{\overline{f_{\geq k_{0}t+2}(z,w)} \right\}\geq k_{0}t+2
,\quad \mathrm{wt}\left\{f_{\rm{nor}}^{\left(k_{0}t+1\right)}(z,w)\right\},\mathrm{wt}\left\{\overline{f_{\rm{nor}}^{\left(k_{0}t+1\right)}(z,w)}\right\}\geq ts+1, \label{1e1se}\end{equation}
where $w$ satisfies (\ref{ec4}).   Furthermore, by (\ref{1e1se}) it follows that
\begin{equation}\begin{split}& \quad\quad\quad\quad\quad\quad\quad\quad\quad\quad\quad\quad\quad P\left(z+\displaystyle\sum_{l=0}^{\left(s-k_{0}\right)\left(t+1\right)}f_{\rm{nor}}^{\left(k_{0}t+l+1\right)}(z,w),\overline{z+\displaystyle\sum_{l=0}
^{\left(s-k_{0}\right)\left(t+1\right)}f_{\rm{nor}}^{\left(k_{0}t+l+1\right)}(z,w)}\right)=P(z,\overline{z})\\&\quad\quad+2\Re\left\{P_{z}(z,\overline{z})\displaystyle\sum_{l=0}^{\left(s-k_{0}\right)\left(t+1\right)}f_{\rm{nor}}^{\left(k_{0}t+l+1\right)}(z,w)\right\} +
2\Re\left\{\displaystyle\sum_{\gamma+\beta=2}^{k_{0}}P_{z^{\gamma}\overline{z}^{\beta}}(z,\overline{z})\left(
\displaystyle\sum_{l=0}^{\left(s-k_{0}\right)\left(t+1\right)}f_{\rm{nor}}^{\left(k_{0}t+l+1\right)}(z,w)\right)^{\gamma}\right. \\& \left.\left(\overline{
\displaystyle\sum_{l=0}^{\left(s-k_{0}\right)\left(t+1\right)}f_{\rm{nor}}^{\left(k_{0}t+l+1\right)}(z,w)}\right)^{\beta}\right\} =P(z,\overline{z})+2\Re\left\{P_{z}(z,\overline{z})f_{\rm{nor}}^{\left(k_{0}t+1\right)}(z,w)\right\}+2\Re\left\{\left(P_{z}(z,\overline{z})+\Theta_{s}^{k_{0}}(z,\overline{z})\right)f_{\geq k_{0}t+2}(z,w)\right\}
\\&\quad\quad\quad\quad\quad\quad\quad\quad\quad\quad\quad\quad\quad\quad\quad \quad\quad\quad\quad\quad\quad\quad\quad\quad\quad\quad\quad\quad\quad\quad \quad\quad\quad\quad\quad\quad\quad\quad\quad\quad\quad\quad\quad\quad\quad \quad   +\Theta_{ts+s+1}^{k_{0}t+k_{0}+1}(z,\overline{z}),\end{split}\end{equation}
where $w$ satisfies (\ref{ecg4}) and $\rm{wt}\left\{\overline{\Theta_{ts+s+1}^{k_{0}t+k_{0}+1}(z,\overline{z})}\right\}\geq ts+s+1$.
 
 In order to track the action of the free parameter $a$, we   prove  

 \bl\label{L3.3se} For all  $m,n\geq 1$ with $m+n\geq k_{0}+1$ and $w$ satisfying (\ref{ec4}), we
have  
\begin{equation} \left(z+f(z,w)\right)^{m}\left(\overline{z+f(z,w)}\right)^{n}=z^{m}\overline{z}^{n}
+2\Re\left\{\Theta_{s}^{k_{0}}(z,\overline{z})
 f_{\geq
k_{0}t+2}(z,w)\right\}+\Theta_{ts+s+1}^{k_{0}t+k_{0}+1}(z,\overline{z}),\end{equation}
where
$\mathrm{wt}\left\{\overline{\Theta_{ts+s+1}^{k_{0}t+k_{0}+1}(z,\overline{z})}\right\}\geq
 ts+s+1$.
 \el
\begin{proof} By the Taylor expansion, we obtain  
\begin{equation}\begin{split}\left(z+f(z,w)\right)^{m}\left(\overline{z+f(z,w)}\right)^{n}&=z^{m}\overline{z}^{n}
+2\Re\left\{\Theta_{s}^{k_{0}}(z,\overline{z})\left(f_{\rm{nor}}^{\left(k_{0}t+1\right)}(z,w)+
 f_{\geq
k_{0}t+2}(z,w)\right)\right\} \\&=z^{m}\overline{z}^{n}
+2\Re\left\{\Theta_{s}^{k_{0}}(z,\overline{z})
 f_{\geq
k_{0}t+2}(z,w)\right\}+\Theta_{ts+s+1}^{k_{0}t+k_{0}+1}(z,\overline{z}),\end{split}\end{equation}
where
$\mathrm{wt}\left\{\overline{\Theta_{ts+s+1}^{k_{0}t+k_{0}+1}(z,\overline{z})}\right\}\geq
 ts+s+1$. \end{proof}

Making computations as in the first studied case,  we obtain
 \bl\label{L091} For $w$ satisfying (\ref{ecg4}) and for all $k>s$, we have  
 \begin{equation} \left(z+f(z,w)\right)^{k}=z^{k}+2\Re \left\{
\Theta_{s}^{k_{0}}(z,\overline{z})f_{\geq
k_{0}t+2}(z,w)\right\}+\Theta_{ts+s+1}^{k_{0}t+k_{0}+1}(z,\overline{z}),\end{equation}
where
$\mathrm{wt}\left\{\overline{\Theta_{ts+s+1}^{k_{0}t+k_{0}+1}(z,\overline{z})}\right\}\geq
 ts+s+1$.
 \el
  \bl\label{L091} For $w$ satisfying (\ref{ecg4}), we have  
\begin{equation} g_{\rm{nor}}^{\left(k_{0}t+k_{0}\right)}(z,w)-2\Re\left\{P\left(f_{\rm{nor}}^{\left(k_{0}t+1\right)}(z,w),\overline{z}\right)\right\}=2\Re \left\{  \left( \alpha   a + \overline{\alpha}\overline{z} \right)  z^{s}w^{t}-aR(z,\overline{z}) w^{t} \right\}+\Theta_{ts+s+1}^{k_{0}t+k_{0}+1}(z,\overline{z}),\end{equation}
where
$\mathrm{wt}\left\{\overline{\Theta_{ts+s+1}^{k_{0}t+k_{0}+1}(z,\overline{z})}\right\}\geq
 ts+s+1$.
 \el
 \begin{proof} It follows easily by the fact that the normalization conditions   (\ref{cn}) remain preserved and by using  the decomposition   (\ref{alpha}).
 \end{proof}
 
 By Lemmas \ref{L3.3se}-\ref{L091} and using (\ref{ecg4se}), in order to track the action of the parameter $a$,    we prove  
 \bl \label{L3.7se}Let $N'_{s}:=ts+s+1$. For all $0\leq j\leq t$ and
$p\in\left[k_{0}t+j\left(s-k_{0}\right)+k_{0},k_{0}t+(j+1)\left(s-k_{0}\right)+k_{0}\right]$, we have  
\begin{equation}\begin{split} g_{\geq p}(z,w)= & 2\Re\left\{z^{sj+s}\left( \left(\left(s-\alpha \right) a - \overline{\alpha}{\overline{a}} \right)\left(1-\frac{s}{\alpha}\right)^{j} 
w^{t-j}+  R(z,\overline{z})\sum_{\beta+\gamma=t-j-1}v_{\beta\gamma}w^{\beta}\left(P(z,\overline{z})\right)^{\gamma}\right)     \right\}\\&+2\Re\left\{\left( P_{z}(z,\overline{z})+
sz^{s-1}+\Theta_{s}^{k_{0}}(z,\overline{z})\right)
f_{\geq p-k_{0}+1}(z,w) \right\}    +a'_{\geq p}(z,\overline{z})-a_{\geq
p}(z,\overline{z})+\Theta_{N'_{s}}^{p}(z,\overline{z}),\end{split}
 \label{30se}
\end{equation}
where
$\mathrm{wt}\left\{\overline{\Theta_{N'_{s}}^{p}(z,\overline{z})}\right\}\geq
 N'_{s}$, $w$ satisfies (\ref{ec4}) and $\alpha$ is defined by (\ref{alpha}).
  \el
\begin{proof} The proof of the lemma has the same structure as the proof of Lemma \ref{L3.7}.

$\bf{Step\hspace{0.1 cm}1.}$ It is very similarly to the first step of the proof of Lemma \ref{L3.7}.

$\bf{Step\hspace{0.1 cm}2.}$ Assume that
we have proved Lemma \ref{L3.7se} for
$p\in\left[k_{0}t+j\left(s-k_{0}\right)+k_{0},k_{0}t+(j+1)\left(s-k_{0}\right)+k_{0}\right]$, for 
$j\in[0,t-1]$. We  prove  Lemma \ref{L3.7se} for
$p\in\left[k_{0}t+(j+1)\left(s-k_{0}\right)+k_{0},k_{0}t+(j+2)\left(s-k_{0}\right)+k_{0}\right]$. Adapting Huang-Yin's strategy\cite{HY2}, we define  
\begin{equation}\begin{split}& f_{\rm{nor}}^{\left(\Lambda\right)}(z,w)=f_{1}^{\left(\Lambda \right)}(z,w)+f_{2}^{\left(\Lambda \right)}(z,w),\\& f_{1}^{\left(\Lambda\right)}(z,w)=-\frac{a}{\alpha}\left(\left(s-\alpha \right) a - \overline{\alpha}\overline{z} \right)\left(1-\frac{s}{\alpha}\right)^{j} z^{sj+s+1}
w^{t-j-1},\quad\mbox{for $\Lambda=k_{0}t+(j+1)\left(s-k_{0}\right)+1$.} \end{split}  \label{35se}\end{equation}

Substituting
$f_{\geq\Lambda}(z,w)=f_{\rm{nor}}^{\left(\Lambda\right)}(z,w)+f_{\geq\Lambda+1}(z,w)$
and
$g_{\geq\Lambda+k_{0}-1}(z,w)=g_{\rm{nor}}^{\left(\Lambda+k_{0}-1\right)}(z,w)+g_{\geq\Lambda+k_{0}}(z,w)$
in (\ref{30se}), we obtain
\begin{equation} \begin{split}  g_{\rm{nor}}^{\left(\Lambda+k_{0}-1\right)}(z,w)+g_{\geq\Lambda+k_{0}}(z,w)=& 2\Re\left\{z^{sj+s}\left(\left(\left(s-\alpha \right) a - \overline{\alpha}\overline{z} \right)\left(1-\frac{s}{\alpha}\right)^{j} 
w^{t-j}+  R(z,\overline{z})\sum_{\beta+\gamma=t-j-1}v_{\beta\gamma}w^{\beta}\left(P(z,\overline{z})\right)^{\gamma}\right)  \right\}\\&+2\Re\left\{
\left(P_{z}(z,\overline{z})+ sz^{s-1}
+\Theta_{s}^{k_{0}}(z,\overline{z})\right)\left(f_{\rm{nor}}^{\left(\Lambda \right)}(z,w)+f_{\geq
\Lambda+1}(z,w)\right)\right\}\\&
 +a'_{\geq\Lambda+k_{0}-1}(z,\overline{z})-a_{\geq\Lambda+k_{0}-1}(z,\overline{z})+\mathbb{P}_{N'_{s}}^{\Lambda+k_{0}-1}(z,\overline{z})
+\Theta_{N'_{s}}^{\Lambda+k_{0}}(z,\overline{z}).\end{split}  \label{38se}\end{equation}

Rewriting  (\ref{38se}) using (\ref{35se}),  we obtain 
\begin{equation}\begin{split}  g_{\rm{nor}}^{\left(\Lambda+k_{0}-1\right)}(z,w)+g_{\geq\Lambda+k_{0}}(z,w)=&
2\Re\left\{\left(P_{z}(z, \overline{z})+
sz^{s-1}+\Theta_{s}^{k_{0}}(z,\overline{z})\right)\left(f_{\geq\Lambda+1 }(z,w)+f^{\left(\Lambda \right)}_{2}(z,w)\right)\right\}
\\&+a'_{\geq\Lambda+k_{0}-1}(z,\overline{z})-a_{\geq\Lambda+k_{0}-1}(z,\overline{z})
+\Theta_{N'_{s}}^{\Lambda+k_{0}}(z,\overline{z})+J'(z,\overline{z}),\end{split}  \label{39}\end{equation}
where we have used the  notation
\begin{equation*}\begin{split}  J'(z,\overline{z})=&2\Re\left \{P_{z}(z,\overline{z})\left(f^{\left(\Lambda \right)}_{1}(z,w)-f^{\left(\Lambda\right)}_{1}\left(z,P(z,\overline{z})\right) \right)\right\}+2\Re \left\{\left(\Theta_{s}^{k_{0}}(z,\overline{z})+sz^{s-1}\right)f_{1}^{\left(\Lambda \right)}(z,w)\right\}\\&  +2\Re\left\{ z^{sj+s}\left(\left(\left(s-\alpha \right) a - \overline{\alpha}\overline{z} \right)\left(1-\frac{s}{\alpha}\right)^{j} 
\left(w^{t-j}-\left(P(z,\overline{z})\right)^{t-j}\right)+  R(z,\overline{z})\sum_{\beta+\gamma=t-j-1}v_{\beta\gamma}\left(w^{\beta}-\left(P(z,\overline{z})\right)^{\beta}\right) \left(P(z,\overline{z})\right)^{\gamma}\right) \right\}
.\end{split}    \end{equation*}

Because the Fischer normalization conditions  (\ref{cn}) remain preserved, it follows that
$$J'(z,\overline{z})=2\Re\left\{z^{s(j+1)+s}\left(\left(\left(s-\alpha \right) a - \overline{\alpha}\overline{z} \right)\left(1-\frac{s}{\alpha}\right)^{j+1} 
w^{t-j-1}+  R(z,\overline{z})\sum_{\beta+\gamma=t-j-2}v_{\beta\gamma}w^{\beta}\left(P(z,\overline{z})\right)^{\gamma}\right)  \right\}+\Theta_{N'_{s}}^{\Lambda+k_{0}}(z,\overline{z}),$$
where $\rm{wt}\left\{\overline{\Theta_{N'_{s}}^{\Lambda+k_{0}}(z,\overline{z})}\right\}\geq N'_{s}$. We    are able to finish the proof, because the multiple of $R(z,\overline{z})$ does not appear when $t=j$ in  (\ref{30se}).
\end{proof}

Collecting the  terms of bidegree $(m,n)$ in $(z,\overline{z})$
with $m+n=ts+s$ and $t=j$ in  (\ref{30se}), we obtain  
\begin{equation}\begin{split} g_{\rm{nor}}^{\left(ts+s\right)}\left(z,P(z,\overline{z})\right)=&2 \Re\left\{\left(\left(s-\alpha \right) a - \overline{\alpha}\overline{z} \right)
\left(1-\frac{s}{\alpha}\right)^{t}z^{ts+s} \right\}+2\Re\left\{P_{z}(z,\overline{z})
f_{\rm{nor}}^{\left(ts+s-k_{0}+1\right)}\left(z,P(z,\overline{z})\right)\right\} \\& 
+a'_{ts+s}(z,\overline{z})-a_{ts+s}(z,\overline{z})+\left(\Theta_{1}\right)_{N'_{s}}^{ts+s}(z,\overline{z}),\end{split} \label{46}\end{equation}
$\rm{wt}\left\{\overline{\left(\Theta_{1}\right)_{N'_{s}}^{ts+s}(z,\overline{z})}\right\}\geq N'_{s}$. By imposing the corresponding Fischer normalization conditions   (\ref{cn}), we find the solution 
$$\left(f_{\rm{nor}}^{\left(ts+s-k_{0}+1\right)}(z,w),\hspace{0.1 cm}g_{\rm{nor}}^{\left(ts+s\right)}(z,w)\right)\hspace{0.1 cm}\mbox{ for
(\ref{46}).}$$ 

 The parameter
$a$ is computed by imposing the corresponding  normalization condition from
(\ref{extra}) in this case.  By composing the map that sends $M$ into
(\ref{ec4}) with the map (\ref{t1se}), we obtain our formal
transformation that sends $M$ into $M'$ up  to degree $ts+s$. 

The uniqueness of the formal transformation (\ref{2800}) can be proven following the lines from \cite{V1}. It is enough to show that any formal equivalence between two real formal surfaces satisfying the Fischer normalization conditions   (\ref{cn}) and (\ref{extra}) is just the identity. The proof follows as in \cite{V1} applying the arguments used during the first step and during the second step of the above construction.

\section{Open Problem}

It would be interesting to construct an example of two real-analytic submanifolds in somplex space that are formally
biholomorphically equivalent, but not biholomorphically equivalent. In $\mathbb{C}^{2}$ we mention that Gong \cite{G3} in the CR singular case, and respectively
Kossovskiy-Shafikov \cite{KS}  in the CR-case, constructed examples of real-analytic submanifolds in  complex space satisfying such property.

\section{Ackowlodgements} This project has been started by me when I was  Ph.D. student at  the School of Mathematics, Trinity College
Dublin, Ireland. I am grateful to Prof. Dmitri Zaitsev for constant support and to Prof. Xiaojun Huang for useful discussions regarding the
Generalization\cite{HY1} of the Theorem of Moser\cite{Mos}.    This project was supported partially by CAPES at the Federal University of Santa Catarina, Brazil. 
The reference \cite{V1} has been fully supported by  Science Foundation Ireland, Grant 06/RFP/MAT 018.

\end{document}